\documentclass[11pt,leqno,a4paper]{article}

\usepackage{amsmath, amssymb, amsthm, amsfonts}
\usepackage{mathrsfs}

\theoremstyle{plain}
\newtheorem{thm}{Theorem}

\newtheorem{lem}{Lemma}

\theoremstyle{remark}
\newtheorem{rem}{Remark}

\numberwithin{equation}{section}

\begin{document}

\title{Inequalities for Euler-Mascheroni constant}
\author{ Hongmin Xu and Xu You }
\date{Jun 26, 2014}

\maketitle



\begin{abstract}
The aim of this paper is to establish new inequalities for the Euler-Mascheroni by the continued fraction method.

\end{abstract}

\section{Introduction}
Euler's Constant was first introduced by Leonhard Euler (1707-1783) in 1734 as the limit of the sequence
\begin{align}
\gamma(n):=\sum_{m=1}^{n}\frac 1m -\ln n.
\end{align}
It is also known as the Euler-Mascheroni constant. There are many famous unsolved problems about the nature of this constant(See e.g. the survey papers or books of R.P. Brent and P. Zimmermann\cite{BZ}, Dence and Dence\cite{DD}, Havil\cite{Ha} and Lagarias\cite{Lag}). For example, it is a long-standing open problem if it is a rational number. A good part of its mystery comes from the fact that the known algorithms converging to $\gamma$ are not very fast, at least, when they are compared to similar algorithms for $\pi$ and $e$.

The sequence $\left(r(n)\right)_{n\in \mathbb{N}}$ converges very slowly toward $\gamma$, like $(2n)^{-1}$. Up to now, many authors are
preoccupied to improve its rate of convergence, see e.g. \cite{CM,DD,De,GI,GS,Lu,Lu1,Mor1,MC} and references therein. We list some main results as follows:
\begin{align*}
\sum_{m=1}^{n}\frac 1m -\ln \left(n+\frac 12\right)=\gamma+O(n^{-2}),
\quad \mbox{(DeTemple, \cite{De})}
\end{align*}

\begin{align*}
\sum_{m=1}^{n}\frac 1m -\ln
\frac{n^3+\frac 32n^2+\frac{227}{240}+\frac{107}{480}}
{n^2+n+\frac{97}{240}}=\gamma+O(n^{-6}), \quad \mbox{(Mortici, \cite{Mor1})}
\end{align*}

\begin{align*}
\sum_{m=1}^{n}\frac 1m -\ln\left(
1+\frac {1}{2n}+\frac{1}{24n^2}-\frac{1}{48n^3}
+\frac{23}{5760n^4}\right)=\gamma+O(n^{-5}), \quad \mbox{(Chen and Mortici, \cite{CM})}
\end{align*}

Recently, Mortici and Chen\cite{MC} provided a very interesting sequence
\begin{align*}
\nu(n)=&\sum_{m=1}^{n}\frac 1m -\frac 12\ln\left(n^2+n+\frac 13\right)\\
&-\left(\frac {-\frac{1}{180}}{\left(n^2+n+\frac 13\right)^2}
+\frac {\frac{8}{2835}}{\left(n^2+n+\frac 13\right)^3}+
+\frac {\frac{5}{1512}}{\left(n^2+n+\frac 13\right)^4}
+\frac {\frac{592}{93555}}{\left(n^2+n+\frac 13\right)^5}
\right),
\end{align*}

and proved
\begin{align}
\lim_{n\rightarrow \infty}n^{12}\left(\nu(n)-\gamma\right)=-\frac{796801}{43783740}.
\end{align}
Hence the rate of the convergence of the sequence $\left(\nu(n)\right)_{n\in \mathbb{N}}$ is $n^{-12}$.

Very recently, by inserting the continued fraction term in (1.1), Lu\cite{Lu} introduced  a class of sequences $\left(r_k(n)\right)_{n\in \mathbb{N}}$(see Theorem 1 below), and showed
\begin{align}
&\frac{1}{72(n+1)^{3}}<\gamma -r_{2}(n)<\frac{1}{72n^{3}},\\
&\frac{1}{120(n+1)^{4}}< r_{3}(n)-\gamma<\frac{1}{120(n-1)^{4}}.
\end{align}
In fact, Lu\cite{Lu} also found $a_4$ without proof. In general, the continued fraction method could provide a better approximation
 than others, and has less numerical computations.

First, we will prove
\begin{thm} For Euler-Mascheroni constant, we have the following convergent sequence
\begin{align*}
r(n)=1+\frac 12 +\cdots+\frac 1n-\ln n-\frac{a_1}{n+\frac{a_2 n}{n+\frac{a_3n}{n+\ddots}}},
\end{align*}
where
$(a_1,a_2,a_4,a_6,a_8,a_{10},a_{12})=\left(
\frac{1}{2},\frac{1}{6},\frac{3}{5},
\frac{79}{126},\frac{7230}{6241},
\frac{4146631}{3833346},
\frac{306232774533}{179081182865}\right)$,
and $a_{2k+1}=-a_{2k}$ for $1\le k\le 6$.

Let
\begin{align*}
R_k(n):=\frac{a_1}{n+\frac{a_2 n}{n+\frac{a_3n}{n+\frac{a_4n}{\frac{\ddots}{n+a_k}}}}},
\end{align*}
(See Appendix for their simple expressions) and
\begin{align*}
r_k(n):=\sum_{m=1}^{n}\frac 1m -\ln n-R_k(n).
\end{align*}
For $1\le k\le 13$, we have
\begin{align}
\lim_{n\rightarrow\infty}n^{k+1}\left(r_k(n)-\gamma\right)=C_k,
\end{align}
where
$(C_1,\cdots,C_{13})=\left(-\frac{1}{12},-\frac{1}{72},\frac{1}{120},
\frac{1}{200},-\frac{79}{25200},-\frac{6241}{3175200},\frac{241}{105840},
\frac{58081}{22018248},-\frac{262445}{91974960},\right.$

$\left. -\frac{2755095121}{892586949408},\frac{20169451}{3821257440},
\frac{406806753641401}{45071152103463200},
-\frac{71521421431}{5152068292800}\right)$.
\end{thm}
\bigskip

\noindent{\bf Open problem}  For every $k\ge 1$,  we have $a_{2k+1}=-a_{2k}$ .
\bigskip

The main aim of this paper is to improve (1.3) and (1.4). We establish the following more precise inequalities.
\begin{thm} Let $r_{10}(n), r_{11}(n)$,  $C_{10}$ and $C_{11}$ be defined in Theorem 1, then
\begin{align}
&C_{10}\frac{1}{(n+1)^{11}}<\gamma -r_{10}(n)<C_{10}\frac{1}{n^{11}},\\
&C_{11}\frac{1}{(n+1)^{12}}< r_{11}(n)-\gamma<C_{11}\frac{1}{n^{12}}.
\end{align}
\end{thm}

\bigskip
\begin{rem} In fact, Theorem 2 implies that
$r_{10}(n)$ is a strictly increasing function of $n$, whereas
$r_{11}(n)$ is a strictly decreasing function of $n$. Certainly, it has the similar inequalities for $r_k(n) (1\le k\le 9)$, we leave these for readers to verify. It is also should be noted that (1.4) cannot deduce the monotony of $r_3(n)$.
\end{rem}
\bigskip
\begin{rem}It is
 worth to pointing out that Theorem 2 provides sharp bounds for harmonic sequence, which are superior to Theorem 3 and 4 of Mortici and Chen\cite{MC}.
\end{rem}

\section{The Proof of Theorem 1}

The following lemma gives a method for measuring the rate of convergence, This Lemma was first used by Mortici\cite{Mor2,Mor3} for constructing asymptotic expansions, or to accelerate some convergences. For proof and other details, see, e.g., \cite{Mor3}.
\begin{lem}
If the sequence $(x_n)_{n\in \mathbb{N}}$ is convergent to zero and there exists the limit
\begin{align}
\lim_{n\rightarrow +\infty}n^s(x_n-x_{n+1})=l\in [-\infty,+\infty]
\end{align}
with $s>1$, then there exists the limit:
\begin{align}
\lim_{n\rightarrow +\infty}n^{s-1}x_n=\frac{l}{s-1}.
\end{align}
\end{lem}

In the sequel, we always assume $n\ge 2$.

We need to find the value $a_1\in \mathbb{R}$ which produces the most accurate approximation of the form
\begin{align}
r_1(n)=\sum_{m=1}^{n}\frac 1m -\ln n-\frac{a_1}{n},
\end{align}
here we note $R_1(n)=\frac{a_1}{n}$. To measure the accuracy of this approximation, we usually say that an approximation (2.3) is better as $r_1(n) -\gamma$ faster converges to zero. Clearly
\begin{align}
r_1(n)-r_{1}(n+1)=\ln \left(1+\frac 1n\right)-\frac {1}{n+1}+\frac{a_1}{n+1}-\frac{a_1}{n}.
\end{align}
It is well-known that for $|x|<1$,
\begin{align*}
\ln(1+x)=\sum_{m=1}^{\infty}(-1)^{m-1}\frac{x^{m}}{m}
\quad \mbox{and}\quad \frac{1}{1-x}=\sum_{m=0}^{\infty}x^m.
\end{align*}
Developing the expression (2.4) into power series expansion in $\frac 1n$, we easily obtain
\begin{align}
r_1(n)-r_{1}(n+1)=\left(\frac 12 -a_1\right)\frac{1}{n^2}+\left(a_1-\frac 23\right)\frac{1}{n^3}+
\left(\frac 34-a_1\right)\frac{1}{n^4}+O\left(\frac{1}{n^5}\right).
\end{align}

From Lemma 1, we see that the rate of convergence of the sequence $\left(r_1(n)-\gamma\right)_{n\in \mathbb{N}}$ is even higher as the value $s$ satisfying (2.1). By lemma 1, we have

(i) If $a_1\neq \frac 12$, then the rate of convergence of the
$\left(r_1(n)-\gamma\right)_{n\in \mathbb{N}}$ is $n^{-1}$, since
\begin{align*}
\lim_{n\rightarrow\infty}n\left(r_1(n)-\gamma\right)=\frac 12-a_1\neq 0.
\end{align*}

(ii) If $a_1= \frac 12$, from (2.5) we have
\begin{align*}
r_1(n)-r_{1}(n+1)=-\frac 16\frac{1}{n^3}+O\left(\frac{1}{n^4}\right).
\end{align*}

Hence the rate  of convergence of the
$\left(r_1(n)-\gamma\right)_{n\in \mathbb{N}}$ is $n^{-2}$, since
\begin{align*}
\lim_{n\rightarrow\infty}n^2\left(r_1(n)-\gamma\right)=-\frac {1}{12}.
\end{align*}

We also observe that the fastest possible sequence $\left(r_1(n)\right)_{n\in \mathbb{N}}$ is obtained only for $a_1=\frac 12$.

Just as Lu\cite{Lu} did, we may repeat the above approach to determine $a_1$ to $a_4$ step by step. However, the computations
become very difficult when $k\ge 5$. In this paper we will use the \emph{Mathematica} software to manipulate symbolic computations.

Let
\begin{align}
r_k(n)=\sum_{m=1}^{n}\frac 1m -\ln n-R_k(n),
\end{align}
then
\begin{align}
r_k(n)-r_k(n+1)=\ln \left(1+\frac 1n\right)-\frac {1}{n+1}+R_k(n+1)-R_k(n).
\end{align}

It is easy to get the following power series
\begin{align}
\ln \left(1+\frac 1n\right)-\frac {1}{n+1}=\sum_{m=2}^{\infty}
(-1)^m\frac{m-1}{m}\frac{1}{n^m}.
\end{align}
Hence the key step is to expand  $R_k(n+1)-R_k(n)$ into power series in $\frac 1n$. Here we use  some examples to explain our method.

{\bf Step 1: }. For example, given $a_1$ to $a_7$, find $a_8$. Define
\begin{align}
R_8(n)=&\frac{\frac 12}{n +\frac{ \frac n6}{
 n + \frac{-\frac n6}{
  n + \frac{\frac{3}{5}*n}{
   n + \frac{\frac{-3}{5}*n}{n + \frac{\frac{79}{126}*n}{n + \frac{\frac{-79}{126}*n}{n + a_8}}}}}}}
\\
=&\frac{-237 + 1405 a_8 + 1800 n + 1740 a_8 n - 630 n^2 + 3780 a_8 n^2 + 3780 n^3}{6 (79 a_8 + 600 a_8 n + 600 n^2 + 790 a_8 n^2 + 1260 a_8 n^3 + 1260 n^4)}.\nonumber
\end{align}

By using the \emph{Mathematica} software(The \emph{Mathematica Program} is very similar to one given in Remark 3 below, however it has a parameter $a_8$), we obtain
\begin{align}
&R_8(n+1)-R_8(n)\\
=&-\frac{ 1}{2 n^2}+ \frac{2}{3 n^3}- \frac{3}{
 4 n^4} + \frac{4}{5 n^5} - \frac{5}{6 n^6}+\frac{6}{7 n^7}-\frac{7}{8 n^8}\nonumber\\
 &+ \frac{360030 - 6241 a_8}{
 396900 n^9}+\frac{-346440 + 24964 a_8 + 6241 a_8^2}{352800 n^{10}} +O\left(\frac{1}{n^{11}}\right).\nonumber
\end{align}

Substituting (2.8) and (2.10) into (2.7),  we get
\begin{align}
r_8(n)-r_8(n+1)=&\left(-\frac 89+\frac{360030 - 6241 a_8}{
 396900 }\right)\frac{1}{n^9}\\
 &+\left(\frac{9}{10}+\frac{-346440 + 24964 a_8 + 6241 a_8^2}{352800 }\right)\frac{1}{n^{10}} +O\left(\frac{1}{n^{11}}\right).\nonumber
\end{align}
The fastest possible sequence $\left(r_8(n)\right)_{n\in \mathbb{N}}$ is obtained only for $a_8=\frac{7230}{6241}$. At the same time, it follows from (2.11),
\begin{align}
r_8(n)-r_8(n+1)=\frac{58081}{2446472}
\frac{1}{n^{10}} +O\left(\frac{1}{n^{11}}\right),
\end{align}
the rate  of convergence of the
$\left(r_8(n)-\gamma\right)_{n\in \mathbb{N}}$ is $n^{-9}$, since
\begin{align*}
\lim_{n\rightarrow\infty}n^9\left(r_8(n)-\gamma\right)=-\frac {58081}{22018248}.
\end{align*}

We can use the above approach to find $a_k(3\le k\le 8)$. Unfortunately, it does not work well for $a_9$. Since $a_3=-a_2$, $a_5=-a_4$ and $a_7=-a_6$. So we may conjecture $a_9=-a_8$. Now let's  check it carefully.

{\bf Step 2: } Check $a_9=-\frac{7230}{6241}$ to $a_{13}=-\frac{306232774533}{179081182865}$.

Let $a_1,\cdots, a_9$, and $R_9(n)$ be defined in Theorem 1. Applying the \emph{Mathematica} software, we obtain
\begin{align}
&R_9(n+1)-R_9(n)\\
=&-\frac{ 1}{2 n^2}+ \frac{2}{3 n^3}- \frac{3}{
 4 n^4} + \frac{4}{5 n^5} - \frac{5}{6 n^6}+\frac{6}{7 n^7}-\frac{7}{8 n^8}+\frac{ 8}{9}\frac{1}{ n^{9}}\nonumber\\
&-\frac{ 9}{10}\frac{1}{ n^{10}}+\frac{736265}{836136}\frac{1}{ n^{11}}+O\left(\frac{1}{n^{12}}\right),\nonumber
\end{align}
which is the desired result. Substituting (2.8) and (2.13) into (2.7),  we get

\begin{align}
r_9(n)-r_9(n+1)=-\frac{262445}{9197496}
\frac{1}{n^{11}} +O\left(\frac{1}{n^{12}}\right),
\end{align}
the rate  of convergence of the
$\left(r_9(n)-\gamma\right)_{n\in \mathbb{N}}$ is $n^{-10}$, since
\begin{align*}
\lim_{n\rightarrow\infty}n^{10}\left(r_9(n)-\gamma\right)=
-\frac{262445}{91974960}.
\end{align*}

Next, we can use the \mbox{\bf Step 1} to find $a_{10}$, and the \mbox{\bf Step} 2 to check $a_{11}$ and $a_{12}$. It should be noted that
Theorem 2 will provide their another proofs for $a_{10}$ and $a_{11}$. So we omit the details here.

Finally, we check $a_{13}=-\frac{306232774533}{179081182865}$.
\begin{align}
&R_{13}(n+1)-R_{13}(n)\\
=&-\frac{ 1}{2 n^2}+ \frac{2}{3 n^3}- \frac{3}{
 4 n^4} + \frac{4}{5 n^5} - \frac{5}{6 n^6}+\frac{6}{7 n^7}-\frac{7}{8 n^8}+\frac{ 8}{9}\frac{1}{ n^{9}}\nonumber\\
&-\frac{ 9}{10}\frac{1}{ n^{10}}+\frac{10}{11}\frac{1}{ n^{11}}-
\frac{11}{12}\frac{1}{ n^{12}}+\frac{12}{13}\frac{1}{ n^{13}}-
\frac{13}{14}\frac{1}{ n^{14}}
\nonumber\\
&+\frac{1903648586623}{2576034146400}\frac{1}{n^{15}}
+O\left(\frac{1}{n^{16}}\right).\nonumber
\end{align}
Substituting (2.8) and (2.15) into (2.7),  one has
\begin{align}
r_{13}(n)-r_{13}(n+1)=-\frac{500649950017}{2576034146400}\frac{1}{n^{15}}
+O\left(\frac{1}{n^{16}}\right).
\end{align}
Since
\begin{align*}
\lim_{n\rightarrow\infty}n^{14}\left(r_{13}(n)-\gamma\right)=
-\frac{71521421431}{5152068292800},
\end{align*}
thus the rate  of convergence of the
$\left(r_{13}(n)-\gamma\right)_{n\in \mathbb{N}}$ is $n^{-14}$.

This completes the proof of Theorem 1.\qed

\bigskip

\begin{rem}
In fact, if the assertion $a_{13}=-\frac{306232774533}{179081182865}$ holds, then the other values $a_j (1\le j\le 12)$ must be true. The following
\emph{Mathematica program} will generate $R_{13}(n+1)-R_{13}(n)$ into power series in $\frac 1n$ with order 16:

Normal[Series$[(R_{13}[n+1]-R_{13}[n])/. \quad n\rightarrow 1/x, \{x, 0, 16\}]]/.\quad x \rightarrow 1/n $
\end{rem}

\begin{rem}
It is a very interesting question to find $a_{k}$ for $k\ge 14$. However, it seems impossible by the above method.
\end{rem}

\section{The Proof of Theorem 2}
Before we prove the Theorem 2, let us give a simple inequality, which plays an important role of the proof.

\begin{lem} Let $f''(x)$ be a continuous function.  If $f''(x)>0$, then
\begin{align}
\int_{a}^{a+1}f(x) dx > f(a+1/2).
\end{align}
\end{lem}
\proof Let $x_0=a+1/2$. By Taylor's formula, we have
\begin{align*}
\int_{a}^{a+1}f(x) dx=&\int_{a}^{a+1}\left(
f(x_0)+f'(x_0)(x-x_0)+\frac 12 f''(\theta_x)(x-x_0)^2\right)dx\\
>&\int_{a}^{a+1}\left(
f(x_0)+f'(x_0)(x-x_0)\right)dx\\
=&f(a+1/2).
\end{align*}
This completes the proof of Lemma 2.\qed

In the sequel, the notation $P_k(x)$ means a polynomial of degree $k$ in $x$ with all of its non-zero coefficients positive, which may be different at each occurrence.

Let's begin to prove Theorem 2. Note $r_{10}(\infty)=0$, it is easy to see
\begin{align}
\gamma -r_{10}(n)=\sum_{m=n}^{\infty}\left(r_{10}(m+1)-r_{10}(m)\right)
=\sum_{m=n}^{\infty}f(m),
\end{align}
where
$$f(m)=\frac{1}{m+1}-\ln \left(1+\frac 1m\right)-R_{10}(m+1)+R_{10}(m).$$
Let $D_1=\frac{2755095121}{6762022344}$. By using the \emph{Mathematica} software, we have
\begin{align*}
f'(x)+D_1\frac{1}{(x+1)^{13}}=-\frac{P_{19}(x)(x-1)+1619906998377
\cdots
5270931}{33810111720 x (1 + x)^{13}P_{10}^{(1)}(x)P_{10}^{(2)}(x)}<0,
\end{align*}
and
\begin{align*}
f'(x)+D_1\frac{1}{\left(x+\frac 12\right)^{13}}=\frac{P_{22}(x)}{
4226263965 x (1 + x)^2 (1 + 2 x)^{13}P_{10}^{(3)}(x)P_{10}^{(4)}(x)}>0.
\end{align*}

Hence, we get the following inequalities for $x\ge 1$,

\begin{align}
D_1\frac{1}{(x+1)^{13}}<-f'(x)<D_1\frac{1}{(x+\frac 12)^{13}}
\end{align}

Applying $f(\infty)=0$, (3.3) and Lemma 2, we get
\begin{align}
f(m)=& -\int_{m}^{\infty}f'(x) dx \le D_1\int_{m}^{\infty}
\left( x+\frac 12 \right)^{-13} dx\\
=&\frac{D_1}{12}\left(m+\frac 12\right)^{-12}\nonumber
\le \frac{D_1}{12}\int_{m}^{m+1}x^{-12} dx.\nonumber
\end{align}

From (3.1) and (3.4) we obtain
\begin{align}
\gamma -r_{10}(n)\le & \sum_{m=n}^{\infty}\frac{D_1}{12}\int_{m}^{m+1}x^{-12} dx\\
=&\frac{D_1}{12}\int_{n}^{\infty}x^{-12} dx
=\frac{D_1}{132}\frac{1}{n^{11}}.\nonumber
\end{align}

Similarly, we also have
\begin{align*}
f(m)=& -\int_{m}^{\infty}f'(x) dx \ge D_1\int_{m}^{\infty}
\left( x+1\right)^{-13} dx\\
=&\frac{D_1}{12}(m+1)^{-12}
\ge \frac{D_1}{12}\int_{m+1}^{m+2}x^{-12} dx,
\end{align*}
and
\begin{align}
\gamma -r_{10}(n)\ge & \sum_{m=n}^{\infty}\frac{D_1}{12}\int_{m+1}^{m+2}x^{-12} dx\\
=&\frac{D_1}{12}\int_{n+1}^{\infty}x^{-12} dx
=\frac{D_1}{132}\frac{1}{(n+1)^{11}}.\nonumber
\end{align}
Combining (3.5) and (3.6) completes the proof of (1.6).

Note $r_{11}(\infty)=0$, it is easy to deduce
\begin{align}
r_{11}(n)-\gamma=\sum_{m=n}^{\infty}\left(r_{11}(m)-r_{11}(m+1)\right)=\sum_{m=n}^{\infty}g(m),
\end{align}
where $$g(m)=\ln \left(1+\frac 1m\right)-\frac{1}{m+1}-R_{11}(m)+R_{11}(m+1).$$
We write $D_2=\frac{20169451}{24495240}$. By using the \emph{Mathematica} software, we have
\begin{align*}
-g'(x)-D_2\frac{1}{(x+1)^{14}}=\frac{P_{18}(x)}{24495240 x^3 (1 + x)^{14}P_{8}^{(1)}(x)P_{8}^{(2)}(x)}>0
\end{align*}
and
\begin{align*}
-g'(x)-D_2\frac{1}{\left(x+\frac 12\right)^{14}}=-\frac{P_{19}(x)(x-1)+4622005677839353997724676307741
}{6123810 x^3 (1 + x)^3 (1 + 2 x)^{14}
P_{8}^{(3)}(x)P_{8}^{(4)}(x)}<0.
\end{align*}
Hence for $x\ge 1$,
\begin{align}
D_2\frac{1}{(x+1)^{14}}<-g'(x)<D_2\frac{1}{\left(x+\frac 12\right)^{14}}.
\end{align}

Applying $g(\infty)=0$, (3.8) and (3.1), we get

\begin{align}
g(m)=& -\int_{m}^{\infty}g'(x) dx \le D_2\int_{m}^{\infty}
\left(x+\frac 12\right)^{-14} dx\\
=&\frac{D_2}{13}\left(m+\frac 12\right)^{-13}
\le \frac{D_2}{13}\int_{m}^{m+1}x^{-13} dx.\nonumber
\end{align}

It follows from (3.7) and (3.9)
\begin{align}
r_{11}(n)-\gamma \le & \sum_{m=n}^{\infty}\frac{D_2}{13}\int_{m}^{m+1}x^{-13} dx\\
=&\frac{D_2}{13}\int_{n}^{\infty}x^{-13} dx
=\frac{D_2}{156}\frac{1}{n^{12}}.\nonumber
\end{align}
Finally,
\begin{align*}
g(m)=& -\int_{m}^{\infty}g'(x) dx \ge D_2\int_{m}^{\infty}
(x+1)^{-14} dx\\
=&\frac{D_2}{13}(m+1)^{-13}
\ge \frac{D_2}{13}\int_{m+1}^{m+2}x^{-13} dx.
\end{align*}
and
\begin{align}
r_{11}(n)-\gamma \ge & \sum_{m=n}^{\infty}\frac{D_2}{13}\int_{m+1}^{m+2}x^{-13} dx\\
=&\frac{D_2}{13}\int_{n+1}^{\infty}x^{-13} dx
=\frac{D_2}{156}\frac{1}{(n+1)^{12}}.\nonumber
\end{align}
Combining (3.10) and (3.11) completes the proof of (1.7).\qed

\bigskip
\begin{rem} As an example, we give the \emph{Mathematica Program} for the proof of the left-hand side of (3.3):

(i) Together[D[$f[x],\{x,1\}$]$+D_1(x+1)^{13}$];

(ii) Take out the numerator $P[x]$ of the above rational function, then manipulate the program: Apart[$P[x]/(x-1)$].
\end{rem}

\section{Competing Interests}
The authors declare that they have no competing interests.

\section{Authors' Contributions}
Hongmin Xu conceived of the algorithm and helped to draft the manuscript. Xu You carried out the design of the program and drafted the manuscript. All authors read and approved the final manuscript.

\section{Acknowledgements}
The authors thank Prof. Xiaodong Cao for his help.

This research of this paper was supported by the National Natural Science Foundation of China (Grant No.11171344) and the Natural Science Foundation of Beijing (Grant No.1112010).
\bigskip

{\bf Appendix}
For the reader's convenience, we rewrite $R_k(n)$ ($k\le 13$) with minimal denominators as following.
\begin{align*}
R_1(n)=&\frac{1}{2n},\\
R_3(n)=&\frac{1}{2 n}-\frac{1}{12}\frac{1}{n^2},\\
R_5(n)=&\frac{1}{2 n}-\frac{5}{6(1+10 n^2)},\\
R_7(n)=&\frac{1}{2 n}-\frac{79}{1200}\frac{1}{n^2}-\frac{147}{400(10
+21 n^2)},\\
R_9(n)=&\frac{1}{2 n}-\frac{7(871+790 n^2)}{20(241+3990 n^2+3318 n^4)},\\
R_{11}(n)=&\frac{1}{2 n}-\frac{52489}{894348}\frac{1}{n^2}
-\frac{1237227621+584280400 n^2}{4471740(3549+13020 n^2+5302 n^4)},\\
R_{13}(n)=&\frac{1}{2 n}-\frac{39577260671+66288226620 n^2+15762446700 n^4}{
 1260(20169451+434410620 n^2+646328298 n^4+150118540 n^6)}.
\end{align*}

\begin{align*}
R_2(n)=&\frac{3}{6n+1},\\
R_4(n)=&\frac{13 + 30 n}{6 (1 + 6 n + 10 n^2)},\\
R_6(n)=&\frac{5 (281 + 348 n + 756 n^2)}{6 (79 + 600 n + 790 n^2 + 1260 n^3)},\\
R_8(n)=&\frac{964337 + 2646000 n + 2599730 n^2 +
 2621220 n^3}{20 (19039 + 144600 n + 315210 n^2 + 303660 n^3 +
   262122 n^4)},\\
R_{10}(n)=&\frac{7 (108237701 + 208886046 n + 523341290 n^2 + 210464400 n^3 +
     230000760 n^4)}{20 (12649849 + 107768934 n + 209431110 n^2 +
     395365320 n^3 + 174158502 n^4 + 161000532 n^5)},\\
R_{12}(n)=&(3604759235968501 + 11032319618513046 n + 17366281558290420 n^2 +
   19958033982902400 n^3 \\
   &+ 7661417445218460 n^4
    +
   4964130389017800 n^5)/(1260 (1058674313539 + 9019254081474 n\\
   & +
     22801779033180 n^2
      + 33088387754520 n^3 + 33925126033722 n^4 +
     13474242079452 n^5 \\
     &+ 7879572046060 n^6)).\\
\end{align*}

\begin{flushleft}

Hongmin Xu\\
Department of Mathematics and Physics, \\
Beijing Institute of Petro-Chemical Technology,\\
Beijing 102617, P. R. China \\
e-mail: xuhongmin@bipt.edu.cn \\

\bigskip

Xu You \\
1. School of Mathematics and System Science, \\
Beijing University of Aeronautics and Astronautics,\\
Beijing 100191, P. R. China \\
2. Department of Mathematics and Physics, \\
Beijing Institute of Petro-Chemical Technology,\\
Beijing 102617, P. R. China \\
e-mail: youxu@bipt.edu.cn
\end{flushleft}
\end{document}